# INCOMPLETE (EVEN AND ODD) TRIGONOMETRIC SPLINES IN THE PROBLEMS OF CONSTRUCTING APPROXIMATE SOLUTIONS OF SECOND ORDER LINEAR DIFFERENTIAL EQUATIONS

**Denysiuk V.P.**[1], **Doctor of Physical and Mathematical Sciences, Professor,**
**Rybachuk L.V.**[2], **Candidate (Ph. D.) of Physical and Mathematical Sciences, Docent**

[1] National Aviation University, Ukraine
[2] National Technical University of Ukraine "Igor Sikorsky Kyiv Polytechnic Institute", Ukraine
(kvomden@ukr.net, rybachuk.liudmyla@lll.kpi.ua)

**Abstract**

The method of constructing approximate solutions of the first boundary value problem for linear differential equations based on incomplete (even and odd) trigonometric splines is considered. The theoretical positions are illustrated by numerical examples.

**Introduction**

Many problems in science and technology often involve solving linear differential equations. Despite the well-developed theoretical methods for solving such equations, finding their exact solutions is difficult and often impossible. In such cases, the problem of finding approximate solutions to such equations is considered; such solutions are found using numerical methods.

Numerical methods for finding approximate solutions to differential equations can be divided into two classes: numerical and numerical-analytical methods.

Numerical methods allow us to find discrete values of the desired approximate solution on the interval of the differential equation. There are many such methods known (see, for example, [1]), among which we mention the Runge-Kutta, Adams, Euler, Krylov, etc. methods and their modifications. Such methods are quite simple to implement on modern computers and, if applicable, lead to satisfactory results.

When applying numerical-analytical methods, an approximate solution is sought in the form of some sufficiently smooth function that depends on several parameters; these parameters are chosen based on certain conditions. Since we will be dealing with methods of this class in the future, let's consider the use of such methods in more detail.

Suppose we have a linear differential operator of the 2nd order of the form

$$L_2 u = u''(x) + p_1(x)u'(x) + p_2(x)u(x), \qquad (1)$$

where $x \in [a,b]$, $p_1(x), p_2(x)$ are some functions that are generally continuous on $[a,b]$.

Consider the equation

$$L_2 u = f(x). \qquad (2)$$

It is known from the theory of differential equations that the general solution $u(x)$ of this equation contains 2 arbitrary constants. To determine these arbitrary constants, additional 2 conditions are added to equation (2). We will consider the case when these conditions are imposed on the values of the desired solution at points $a$ and $b$ of the segment $[a,b]$, and have the form

$$u(a) = a_0; \quad u(b) = b_0. \qquad (3)$$

Such conditions are called boundary conditions and are said to be the first boundary value problem [2].

Let $C_{[a,b]}^k$ be the class of functions that have a continuous derivative of order $k$ on $[a,b]$. Note that the classes $C_{[a,b]}^0$ contain continuous functions, and the classes $C_{[a,b]}^{-1}$ contain piecewise continuous functions with a finite number of discontinuities of the first kind.

Let's look for a solution to this equation along with the boundary conditions in the form of a function $\varphi(x,\alpha_1,\alpha_2,...,\alpha_N)$ that satisfies the conditions:

a) $\varphi(x,\alpha_1,\alpha_2,...,\alpha_N) \in C_{[a,b]}^r$, where $r \geq 2$;

b) $\varphi(a,\alpha_1,\alpha_2,...,\alpha_N) = a_0;\quad \varphi(b,\alpha_1,\alpha_2,...,\alpha_N) = b_0;$

and $\alpha_1,\alpha_2,...,\alpha_N$ are some unspecified parameters ($N = 1,2,...$).

Substituting the function $\varphi(x,\alpha_0,\alpha_1,...,\alpha_N)$ into (2), we get a residual

$$L_n\varphi(x,\alpha_0,\alpha_1,...,\alpha_N) - f(x) = \varepsilon(x,\alpha_0,\alpha_1,...,\alpha_N).\qquad(4)$$

By minimizing the residual in one sense or another, it is easy to obtain the conditions for determining the parameters of $\alpha_0,\alpha_1,...,\alpha_N$.

There are various methods for minimizing the residual: the collocation method, the least squares method (continuous and discrete variants), the method of moments and its variants, the subdomain method [3], etc. In this paper, without losing generality, we will limit ourselves to the collocation method. As is known, this method consists of the following.

A certain system of points $a \leq x_1 < x_2 < ... < x_N \leq b$ is set on $[a,b]$ and it is require that the following conditions be met

$$\varepsilon(x_i,\alpha_1,\alpha_2,...,\alpha_N) = 0,\quad (i=1,2,...,N).\qquad(5)$$

From this system of $N$ equations, it is possible to determine $N$ parameters $\alpha_0,\alpha_1,...,\alpha_N$.

When applying numerical-analytical methods for constructing approximate solutions, the method of constructing the function $\varphi(x,\alpha_1,\alpha_2,...,\alpha_N)$ plays an important role. In many cases, a generalized polynomial is chosen as this function

$$P_N(x) = \alpha_1\omega_2(x) + \alpha_2\omega_2(x) + ... + \alpha_N\omega_N(x),\qquad(6)$$

where $\omega_1(x)$, $\omega_2(x)$, ..., $\omega_N(x)$ is a system of linearly independent functions.

As a system of functions $\omega_0(x)$, $\omega_1(x)$, ..., $\omega_N(x)$, systems of polynomial functions orthogonal to $[a,b]$ in the continuum or discrete sense, systems of artificially constructed functions, etc. are often chosen; it is clear that these systems of functions have the necessary differential properties. Recently, polynomial B-splines have been chosen as such systems of functions.

The simplest way to use splines in the problem of constructing approximate solutions to differential equations is to use numerical methods to obtain the values of the desired solution at the nodes of certain grids; these values are then interpolated by splines of certain orders. This method looks quite natural and will not be discussed here.

In [4], an algorithm for using polynomial splines to construct approximate solutions to differential equations was probably first presented; fundamental splines were introduced for both the function itself and its derivatives. However, in practice, this algorithm is not widely used, since the construction of such fundamental splines is associated with significant difficulties.

Subsequently, based on the representation of polynomial splines by B-splines, classes of algorithms were developed to construct approximate solutions to many types of differential equations; in most cases, these algorithms use the collocation method in one way or another and are reduced to solving a system of algebraic equations with respect to the parameters that define polynomial splines [5].

In [6], classes of trigonometric splines were introduced. The properties of these splines are similar to those of simple polynomial splines and were studied in [7], [8]. This naturally raises the question of developing methods for constructing approximate solutions to differential equations using such splines. Given the above, the research problem can be formulated as follows.

## Purpose of the work

Development of a method for constructing approximate solutions of the first boundary value problem for linear differential equations based on incomplete (even and odd) trigonometric splines.

## The main part

In [8], three forms of representation of trigonometric splines were given: through discrete (finite) Fourier series, through trigonometric B-splines, and through fundamental trigonometric splines.

When developing methods for constructing approximate solutions to boundary value problems for linear differential equations based on trigonometric splines, the form of representation of trigonometric splines and the interpretation of the parameters of trigonometric splines in each case play an important role. Thus, in the case of representation of trigonometric splines through discrete (finite) Fourier series, discrete Fourier coefficients act as defining parameters; both periodic and non-periodic splines can be constructed.

In the case of representing splines through trigonometric B-splines, the defining parameters are the coefficients at these B-splines, which are difficult to interpret in any reasonable way; at the same time, only periodic splines are currently plotted.

Finally, in the case of representing trigonometric splines through fundamental splines, the determining parameters are the values of the desired solution at the nodes of the given collocation grids; both periodic and non-periodic splines can be constructed.

Although all three forms can be used for methods of constructing approximate solutions, in our opinion, the most convenient form is the one using fundamental trigonometric splines. This is because using the values of the desired solution as parameters greatly simplifies the enforcement of the boundary conditions for this solution. Therefore, in the following, we will limit ourselves to considering the case when the desired solution is constructed using fundamental trigonometric splines.

Let $[a,b]$ ($a > b$) be a certain segment of the real axis, and let there be a certain grid $\Delta_N = \{x_j\}_{j=1}^N$ on this segment, and $a \leq x_1 < x_2 < ... < x_N \leq b$. The system of functions $\varphi_k(x)$ ($k = 1, 2, ..., N$) on the grid $\Delta_N$ is called fundamental if the following conditions are met [9]

$$\varphi_k(x_j) = \begin{cases} 1, & k = j; \\ 0, & k \neq j, \quad k, j = 1, 2, ..., N. \end{cases} \tag{7}$$

Given a system of fundamental functions, the function $\Phi_N(x)$ interpolating a given function $f(x)$ at the nodes of the grid $\Delta_N$ can be represented as

$$\Phi_N(x) = \sum_{k=1}^{N} f(x_k)\varphi_k(x).$$

Similarly, systems of fundamental polynomial splines are introduced [10].

Systems of trigonometric fundamental splines are of three types: complete fundamental splines, as well as even and odd fundamental splines [8].

Complete fundamental splines are periodic functions and can be used to find periodic solutions of linear differential equations if such solutions exist; recall that the existence of such solutions is discussed, for example, in [2]. In this paper, we will not consider the problem of finding periodic solutions of differential equations.

Even and odd fundamental splines can be used to find nonperiodic solutions to differential equations. Let us consider this issue in more detail.

Let us define uniform grids $\Delta 2_N^{(I)} = \{x_j^{(I)}\}_{j=1}^N$ on $[0, \pi]$, where $x_j^{(0)} = \frac{\pi}{N-1}(j-1)$, $x_j^{(1)} = \frac{\pi}{2N}(2j-1)$ and $\Delta 3_N^{(I)} = \{x_j^{(I)}\}_{j=1}^N$, where $x_j^{(0)} = \frac{\pi}{N+1}j$, $x_j^{(1)} = \frac{\pi}{2N}(2j-1)$, $I$ is the grid indicator ($I = 0,1$), and $N = 1, 2, 3, ...$. We will not introduce separate notations for the nodes of these grids, since in each case we will determine which grid we are talking about.

Non-periodic even and odd trigonometric splines are constructed on the grids $\Delta 2_N^{(I)}$ and $\Delta 3_N^{(I)}$, respectively [8]. Even and odd trigonometric splines have different properties; in particular, even trigonometric splines take non-zero values at the points $0$ and $\pi$, and their derivatives take zero values at these points. Odd splines, on the other hand, take on zero values at the points $0$ and $\pi$, and their derivatives at these points take on non-zero values. Taking this into account, and wanting to deal with uniform grids, we will limit ourselves to considering even splines with the same stitching and interpolation grids $\Delta 2_N^{(0)}$ and odd splines with the same stitching and interpolation grids $\Delta 3_N^{(1)}$.

Here are the expressions for constructing fundamental even splines on the grid $\Delta 2_N^{(0)}$; here and further, in the notation of even and odd fundamental splines, we will omit the indices of the stitching and interpolation grids, as well as the index of the type of convergence factors, which is equal to $\sigma 1(r, j) = j^{-(1+r)}$ ($r, j$ are natural numbers).

The paired fundamental splines $stc(r,q,k,t)$ of order $r$ ($r = 1, 2, ...$), and their derivatives of order $q$ ($q = 1, 2, ..., r-1$) on the grid $\Delta 2_N^{(0)}$ are calculated by the expressions

$$stc(r,q,k,t) = $$
$$= \frac{2}{N-1}\left[.5 \cdot I(q) + \sum_{j=1}^{N-2} \frac{c(\sigma 1, r, q, j, t)\cos(x_k^{(0)})}{H(\sigma 1, r, j)} + .5\frac{c(\sigma 1, r, q, N-1, t)\cos(x_k^{(0)})}{H(\sigma 1, r, N-1)}\right], \quad (8)$$

where

$$I(q) = \begin{cases} 1, & \text{if } q = 0, \\ 0 & \text{otherwise} \end{cases},$$

$$H(\sigma 1, r, k) = \sigma 1(r, k) + \sum_{m=1}^{\infty}\left[\sigma 1(r, 2m(N-1) + k) + \sigma 1(r, 2m(N-1) - k)\right],$$

$$c(\sigma 1, r, q, j, t) = \sigma 1(r, j) j^q \cos(jt + q\frac{\pi}{2}) +$$
$$+ \sum_{m=1}^{\infty}\left[\sigma 1(r, 2m(N-1) + j)(2m(N-1) + j)^q \cos((2m(N-1) + j)t + q\frac{\pi}{2}) +\right.$$
$$\left. + \sigma 1(r, 2m(N-1) - j)(2m(N-1) - j)^q \cos((2m(N-1) - j)t + q\frac{\pi}{2})\right].$$

Graphs of paired trigonometric splines are shown in Figure 1 ($r = 1$ and $N = 9$ for clarity). Note that here and in the following graphs, the vertical lines coincide with the nodes of the grid $\Delta 2_N^{(0)}$.

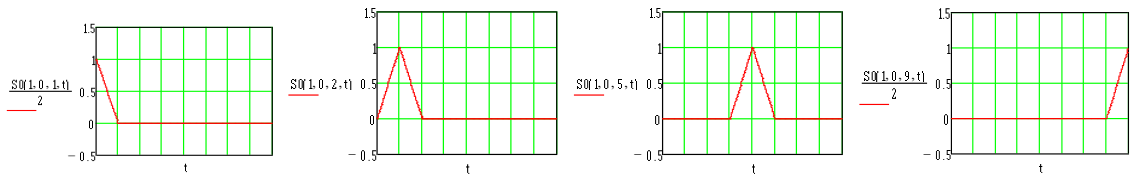

Fig. 1 Paired fundamental splines on the grid $\Delta 2_N^{(0)}$; $r = 1$, $N = 9$

Similarly, the odd fundamental splines $sts^{(0)}(r,q,k,t)$ of order $r$ ($r = 1, 2, ...$), and their derivatives of order $q$ ($q = 1, 2, ..., r-1$) on the grid $\Delta 3_N^{(0)}$, are calculated by the expressions

$$sts^{(0)}(r,q,k,t) = \frac{2}{N+1}\sum_{j=1}^{N}\frac{s(\sigma 1, r, q, j, t)\sin(jx_k^{(0)})}{H(\sigma 1, r, j)}, \quad (9)$$

where

$$H(\sigma 1, r, k) = \sigma 1(r, k) + \sum_{m=1}^{\infty}\left[\sigma 1(r, 2m(N+1) + k) + \sigma 1(r, 2m(N+1) - k)\right]$$

$$s(\sigma 1, r, q, j, t) = \sigma 1(r, j) j^q \sin(jt + q\frac{\pi}{2}) +$$

$$+ \sum_{m=1}^{\infty} (-1)^m \left[ \sigma 1((r, 2m(N+1) + j))((2m(N+1) + j)^q \sin((2m(N+1) + j)t + q\frac{\pi}{2}) - \right.$$

$$\left. - \sigma 1((r, 2m(N+1) - j))((2m(N+1) - j)^q \sin((2m(N+1) - j)t + q\frac{\pi}{2}) \right].$$

As before, we present the graphs of odd trigonometric splines (Figure 2), assuming $r = 1$, $N = 9$. Note that the vertical lines in these graphs coincide with the nodes of the grid $\Delta 3_N^{(0)}$.

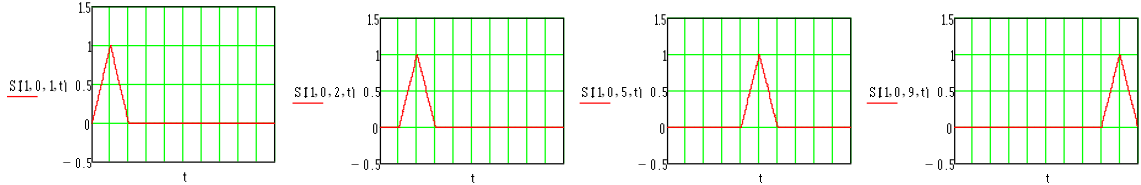

Fig. 2. Odd fundamental splines $sts^{(0)}(1, 0, k, t)$ on the grid $\Delta 3_N^{(0)}$; $r = 1$, $N = 9$

On the grid $\Delta 3_N^{(1)}$, the odd fundamental splines $sts^{(1)}(r, q, k, t)$ of order $r$ ($r = 1, 2, ...$), and their derivatives of order $q$ ($q = 1, 2, ..., r - 1$), are calculated by the expressions

$$sts^{(1)}(r, q, k, t) = \frac{2}{N} \left[ \sum_{j=1}^{N-1} \frac{s(\sigma 1, r, q, j, t) \sin(j x_k^{(1)})}{H(\sigma 1, r, j)} + .5 \frac{s(\sigma 1, r, q, N, t) \sin(N x_k^{(1)})}{H(\sigma 1, r, N)} \right], \qquad (10)$$

where

$$H(\sigma 1, r, k) = \sigma 1(r, k) + \sum_{m=1}^{\infty} \left[ \sigma 1(r, 2mN + k) + \sigma 1(r, 2mN - k) \right];$$

$$s(\sigma 1, r, q, j, t) = \sigma 1(r, j) j^q \sin(jt + q\frac{\pi}{2}) +$$

$$+ \sum_{m=1}^{\infty} (-1)^m \left[ \sigma 1(r, 2mN + j)(2mN + j)^q \sin((2mN + j)t + q\frac{\pi}{2}) - \right.$$

$$\left. - \sigma 1((r, 2mN - j)(2mN - j)^q \sin((2mN - j)t + q\frac{\pi}{2}) \right].$$

Graphs of odd trigonometric splines are shown in Figure 3, assuming $r = 1$, $N = 9$. Note that the vertical lines in these graphs coincide with the nodes of the grid $\Delta 3_N^{(0)}$.

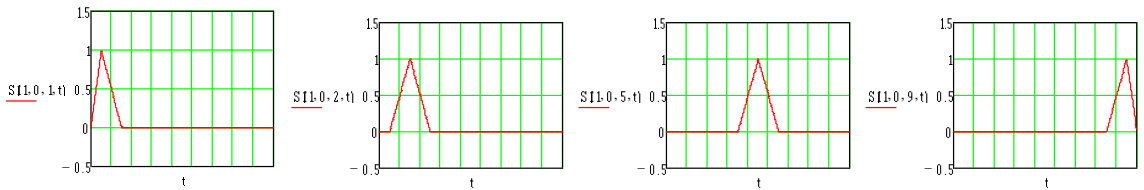

Fig. 3. Odd fundamental splines $sts^{(1)}(1, 0, k, t)$ on the grid $\Delta 3_N^{(0)}$; $r = 1$, $N = 9$

Given the even and odd fundamental splines $stc(r, q, k, t)$, $sts^{(0)}(1, 0, k, t)$ and $sts^{(1)}(1, 0, k, t)$, the splines $Stc(r, q, t)$ and $Sts(r, q, t)$ interpolating a continuous function $f(x)$ at the nodes of the grids $\Delta 2_N^{(0)}$, $\Delta 3_N^{(0)}$ and $\Delta 3_N^{(1)}$ can be represented as

$$Stc(r,q,t) = .5f(x_1)stc(r,q,1,t) + \sum_{k=2}^{N-1} f(x_k)stc(r,q,k,t) + .5f(x_N)stc(r,q,N,t); \qquad (11)$$

$$Sts^{(0)}(r,q,t) = \sum_{k=1}^{N} f(x_k)sts^{(0)}(r,q,k,t), \qquad (12)$$

$$Sts^{(1)}(r,q,t) = \sum_{k=1}^{N} f(x_k)sts^{(1)}(r,q,k,t). \qquad (13)$$

It is clear that interpolation splines $Stc(r,q,k,t)$, $Sts^{(0)}(r,q,k,t)$, $Sts^{(1)}(r,q,k,t)$ and fundamental splines $stc(r,q,k,t)$, $sts^{(0)}(r,q,k,t)$, $sts^{(1)}(r,q,k,t)$ belong to the class $C_{[0,\pi]}^{r-1}$.

Let us now consider the first boundary value problem for equation (2) with boundary conditions (3), which, as we have already mentioned, is formulated as follows. Find the function $u(x)$, $u(x) \in C_{[0,\pi]}^{r}$ ($r \geq 2$), $x \in [0,\pi]$, which satisfies the equation

$$u''(x) + p_1(x)u'(x) + p_2(x)u(x) = f(x) \qquad (14)$$

and edge conditions

$$u(0) = a_0; \quad u(\pi) = b_0. \qquad (15)$$

There are three ways to find an approximate solution $u^*(x)$ to this problem.

a) On the grid $\Delta 2_N^{(0)}$ via paired splines

$$u^*(x) = .5\alpha_1 stc(r,q,1,x) + \sum_{k=2}^{N-1} \alpha_k stc(r,q,k,x) + .5\alpha_N stc(r,q,N,x), \qquad (16)$$

where $r \geq 3$ and the parameters $\alpha_k$ ($k = 2,2,...,N-1$) are determined from the collocation conditions.

It is necessary to take into account the fact that when finding an approximate solution in the form of (16), additional boundary conditions of the form $u'(0) = u'(\pi) = 0$ are also fulfilled. This is due to the fact that when differentiating (16), the functions $\cos(x)$ are transformed into functions $-\sin(x)$, which turn to 0 at the ends of the segment $[0,\pi]$.

b) On the grid $\Delta 3_N^{(0)}$ through odd splines in the form

$$u^{*(0)}(x) = \sum_{k=1}^{N} \alpha_k sts^{(0)}(r,q,k,x), \qquad (17)$$

and on the grid $\Delta 3_N^{(1)}$ in the form of

$$u^{*(1)}(x) = \sum_{k=1}^{N} \alpha_k sts^{(1)}(r,q,k,x), \qquad (18)$$

where $r \geq 3$, and the parameters $\alpha_k$ ($k = 1,2,...,N$) are determined from the collocation conditions. Note that such splines are conveniently used for problems with zero boundary conditions.

In cases where the solution to the equation is sought on the segment $[a,b]$, this segment is mapped to the segment $[0,\pi]$ by linearly replacing the variables.

Consider the case when the solution to problems (14), (15) is sought in the form (17). It is clear that in this case, the nodes $x_1 = 0$ and $x_N = \pi$ of the $\Delta 2_N^{(0)}$ mesh coincide with the ends of the segment $[0,\pi]$, and, therefore, it follows from (15) that $\alpha_1 = a_0$ and $\alpha_N = b_0$. To determine the parameters of $\alpha_2, \alpha_2, ..., \alpha_{N-1}$, we substitute (17) into (14) and equate the values of the left and right parts at the nodes $x_2, x_2, ..., x_{N-1}$ of the grid $\Delta 2_N^{(0)}$; as a result, we obtain the following system of linear algebraic equations with respect to the desired parameters $\alpha_2, \alpha_2, ..., \alpha_{N-1}$:

$$.5a_0 A(r,1,x_j) + \sum_{k=2}^{N-1} \alpha_k A(r,k,x_j) + .5b_0 A(r,N,x_j) = f(x_j) \quad (j = 2,3,...,N-1), \qquad (19)$$

where

$$A(r,k,x_j) = stc(r,2,k,x_j) + p_1(x_j) stc(r,1,k,x_j) + p_2(x_j) stc(r,0,k,x_j)$$
$$(k = 1,2,...,N; \quad j = 2,3,...,N-1).$$

To illustrate these provisions, we will provide examples.

Example 1. Let us consider a boundary value problem

$$u''(x) + \frac{C}{1+x} u'(x) - \frac{x}{1+x} u(x) = \frac{C - 2 - x^2(1+x)}{(1+x)^3}, \quad x \in [0,1];$$
$$u(0) = 0; \quad u(1) = .5, \quad (C = 0, 1, 10).$$

The exact solution to this problem is

$$u0(x) = \frac{x}{1+x}.$$

By composing a system of equations of the type (19) at $r = 3$, and solving it, we obtain approximate solutions $u1(x)$ of this problem for the values $C = 0, 1, 10$. Graphs of exact solution and approximate solutions for different values of the parameter $C$ are shown in Figure 4.

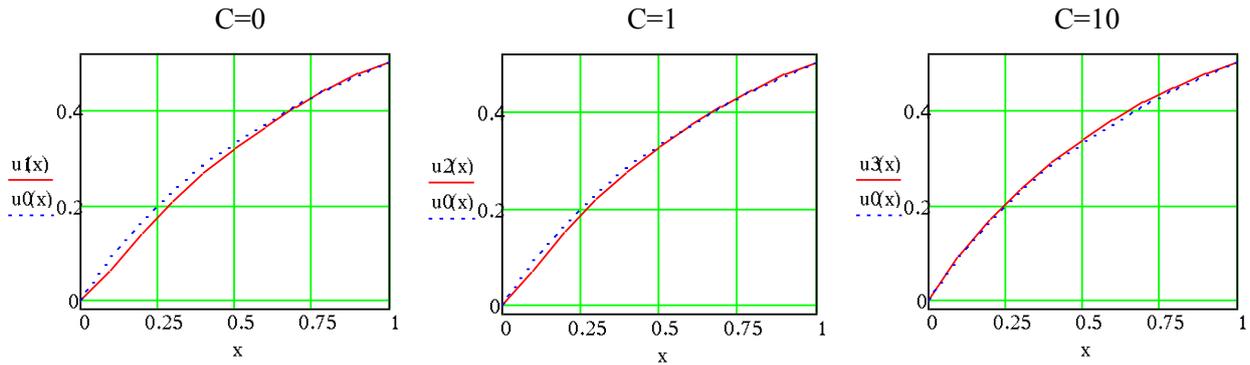

Fig. 4. The exact solution $u0(x)$ and approximate solutions $u1(x)$, $u2(x)$ and $u3(x)$, corresponding to the values of the parameter $C = 0, 1, 10$

Note that the largest deviations of the exact solution from the approximate solution do not exceed 0.048 when C = 0, 0.044 when C = 1, and 0.042 when C = 10.

Example 2. Let us consider a boundary value problem

$$u''(x) + u(x) = \cos(x)\cos(2x), \quad x \in [0,\pi];$$
$$u(0) = 1; \quad u(\pi) = -1.$$

The exact solution to this problem is

$$u0(x) = 1.0625\cos(x) - .4\sin(x) - .0625\cos(3x) + .25x\sin(x).$$

By composing systems of equations of the type (19) at $r = 3, 4, 5$, and solving them, we obtain approximate solutions $u3(x)$, $u4(x)$, $u5(x)$ of this problem corresponding to the given values of the parameter $r$. Graphs of exact solution and approximate solutions for different values of the parameter $r$ are shown in Figure 5.

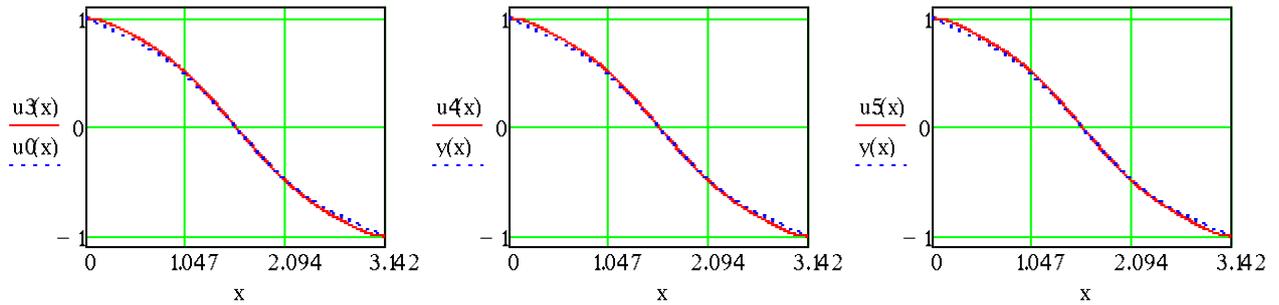

Fig. 5. The exact solution $u0(x)$ and approximate solutions $u3(x)$, $u4(x)$ and $u5(x)$, corresponding to the values of the parameter $r = 3, 4, 5$

The largest deviations of the exact solution from the approximate solution do not exceed 0.048 at $r = 3$, 0.045 at $r = 4$, and 0.043 at $r = 5$.

Let us now consider the case when the solution to problems (14), (15) is sought in the form (18), with (15) having the form

$$u(0) = 0; \quad u(\pi) = 0. \tag{20}$$

It is clear that the nodes $x_1$ and $x_N$ of the $\Delta 3_N^{(0)}$ mesh do not coincide with the ends of the segment $[0, \pi]$. Since (18) is equal to zero at points $0$ and $\pi$, the boundary conditions are automatically satisfied. By drawing up a system of equations of the form

$$\sum_{k=1}^{N} \alpha_k A(r, k, x_j) = f(x_j), \; (j = 2, 3, \ldots, N-1), \tag{21}$$

where

$$A(r, k, x_j) = sts^{(0)}(r, 2, k, x_j) + p_1(x_j) stc^{(0)}(r, 1, k, x_j) + p_2(x_j) stc^{(0)}(r, 0, k, x_j),$$
$$(k = 1, 2, \ldots, N; \quad j = 1, 2, 3, \ldots, N),$$

and finding the value of $\alpha_k$ ($k = 1, 2, \ldots, N$), we construct an approximate solution $u^{*(0)}(x)$.

Reasoning similarly, and using the splines $sts^{(1)}(r, q, k, x)$ to construct the values $A(r, k, x_j)$, it is easy to construct the solution $u^{*(1)}(x)$.

As before, to illustrate the above provisions, we will provide an example.

Example 3. Let us consider a boundary value problem

$$u''(x) + u(x) = -x, \quad x \in [0, \pi] \quad (x \in [0,1]),$$
$$u(0) = 0; \quad u(1) = 0.$$

The exact solution to this problem is

$$u(x) = \frac{\sin(x)}{\sin(1)} - t.$$

By composing systems of equations of the type (21) at $r = 3, 4, 5$, and solving them, we obtain approximate solutions $u3^{(0)}(x)$, $u4^{(0)}(x)$, $u5^{(0)}(x)$ of this problem corresponding to the given values of the parameter $r$. Graph of exact $u(x)$ and approximate $u4^{(0)}(x)$ solutions to this problem is shown in Figure 6; we also show the errors of other solutions for different values of the parameter $r$.

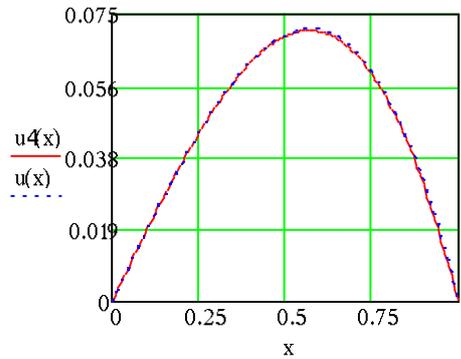

Fig. 6. The exact solution $u(x)$ and the approximate solution $u4^{(0)}(x)$, corresponding to the values of the parameter $r = 4$

The largest deviations of the exact solution from the approximate solution do not exceed 0.00157 at $r = 3$, 0.00105 at $r = 4$, and 0.00098 at $r = 5$.

Similar calculations were performed for the approximate solutions $u3^{(1)}(x)$, $u4^{(1)}(x)$, $u5^{(1)}(x)$; the largest deviations of the exact solution from the approximate solution do not exceed 0.0005 at $r = 3, 4$, but increase to 0.00063 at $r = 5$.

## Conclusions

1. The methods for constructing approximate solutions of the first boundary value problem for linear differential equations based on incomplete (even and odd) trigonometric splines are developed.
2. Paired trigonometric splines should be used in cases where the first derivative of the solution satisfies zero boundary conditions.
3. It is advisable to use odd trigonometric splines for the first boundary value problem with zero boundary conditions.
4. Examples are provided to illustrate the theoretical positions.
5. It is shown that the error of the calculation results depends to a certain extent on the parameter $r$, which characterizes the differential properties of the trigonometric splines used.
6. Undoubtedly, the proposed methods for constructing approximate solutions of the first boundary value problem for linear differential equations based on incomplete (even and odd) trigonometric splines require further research.